\input amstex
\documentstyle{amsppt}
%
\catcode`@=11
\redefine\output@{%
  \def\break{\penalty-\@M}\let\par\endgraf
  \ifodd\pageno\global\hoffset=105pt\else\global\hoffset=8pt\fi  
  \shipout\vbox{%
    \ifplain@
      \let\makeheadline\relax \let\makefootline\relax
    \else
      \iffirstpage@ \global\firstpage@false
        \let\rightheadline\frheadline
        \let\leftheadline\flheadline
      \else
        \ifrunheads@ 
        \else \let\makeheadline\relax
        \fi
      \fi
    \fi
    \makeheadline \pagebody \makefootline}%
  \advancepageno \ifnum\outputpenalty>-\@MM\else\dosupereject\fi
}
\def\Beta{\mathchar"0\hexnumber@\rmfam 42}
\catcode`\@=\active
\nopagenumbers
\chardef\textvolna='176
\def\negskp{\hskip -2pt}

\def\Ker{\operatorname{Ker}}

\def\Sym{\operatorname{Sym}}
\def\compos{\,\raise 1pt\hbox{$\sssize\circ$} \,}

\font\eightrm=cmr8
\def\LT{\operatorname{\text{\eightrm LT}}}

\accentedsymbol\hatgamma{\kern 2pt\hat{\kern -2pt\gamma}}
\accentedsymbol\checkgamma{\kern 2.5pt\check{\kern -2.5pt\gamma}}
\def\blue#1{#1}

\catcode`#=11\def\diez{#}\catcode`#=6
\catcode`&=11\catcode`&=4
\catcode`_=11\def\podcherkivanie{_}\catcode`_=8
\catcode`~=11\def\volna{~}\catcode`~=\active
\def\mycite#1{\cite{\blue{#1}}\immediate\special{ps:
     ShrHPSdict begin /ShrBORDERthickness 0 def}}
\def\myciterange#1#2#3#4{\cite{\blue{#2#3#4}}\immediate\special{ps:
     ShrHPSdict begin /ShrBORDERthickness 0 def}}
\def\mytag#1{%
    \tag#1}
\def\mythetag#1{\thetag{\blue{#1}}\immediate\special{ps:
     ShrHPSdict begin /ShrBORDERthickness 0 def}}
\def\myrefno#1{\no#1}
\def\myhref#1#2{\blue{#2}\immediate\special{ps:
     ShrHPSdict begin /ShrBORDERthickness 0 def}}
\def\myEarXivlink{\myhref{http://arXiv.org}{http:/\negskp/arXiv.org}}

\def\mytheorem#1{\csname proclaim\endcsname{Theorem #1}}
\def\mytheoremwithtitle#1#2{\csname proclaim\endcsname{Theorem #1#2}}
\def\mythetheorem#1{\blue{#1}\immediate\special{ps:
     ShrHPSdict begin /ShrBORDERthickness 0 def}}
\def\mylemma#1{\csname proclaim\endcsname{Lemma #1}}
\def\mylemmawithtitle#1#2{\csname proclaim\endcsname{Lemma #1#2}}
\def\mythelemma#1{\blue{#1}\immediate\special{ps:
     ShrHPSdict begin /ShrBORDERthickness 0 def}}
\def\mycorollary#1{\csname proclaim\endcsname{Corollary #1}}

\def\mydefinition#1{\definition{Definition #1}}
\def\mythedefinition#1{\blue{#1}\immediate\special{ps:
     ShrHPSdict begin /ShrBORDERthickness 0 def}}
\def\myconjecture#1{\csname proclaim\endcsname{Conjecture #1}}
\def\myconjecturewithtitle#1#2{\csname proclaim\endcsname{Conjecture #1#2}}

\font\eightcyr=wncyr8

\pagewidth{360pt}
\pageheight{606pt}
\topmatter
\title
On an ideal of multisymmetric polynomials associated 
with perfect cuboids. 
\endtitle
\rightheadtext{On an ideal of multisymmetric polynomials \dots}
\author
Ruslan Sharipov
\endauthor
\address Bashkir State University, 32 Zaki Validi street, 450074 Ufa, Russia
\endaddress
\email\myhref{mailto:r-sharipov\@mail.ru}{r-sharipov\@mail.ru}
\endemail
\abstract
    A perfect Euler cuboid is a rectangular parallelepiped with integer edges,
with integer face diagonals, and with integer space diagonal as well. Finding 
such parallelepipeds or proving their non-existence is an old unsolved mathematical 
problem. Algebraically the problem is described by a system of Diophantine 
equations. Symmetry approach to the cuboid problem is based on the natural $S_3$
symmetry of its Diophantine equations. Factorizing these equations with respect 
to their $S_3$ symmetry, one gets some certain ideal within the ring of multisymmetric
polynomials. In the present paper this ideal is completely calculated and presented
through its basis. 
\endabstract
\endtopmatter
\TagsOnRight
\document

\head
1. Introduction.
\endhead
     The search for perfect cuboids extends from now back to the year of 1719 
(see\linebreak \myciterange{1}{1}{--}{39}), though one needs only to solve
a very small system of Diophantine equations with respect to seven integer 
variables $x_1$, $x_2$, $x_3$, $d_1$, $d_2$, $d_3$, and $L$:
$$
\xalignat 2
&\hskip -2em
(x_1)^2+(x_2)^2-(d_3)^2=0,
&&(d_3)^2+(x_3)^2-L^2=0,\\
&\hskip -2em
(x_2)^2+(x_3)^2-(d_1)^2=0,
&&(d_1)^2+(x_1)^2-L^2=0,
\mytag{1.1}\\
&\hskip -2em
(x_3)^2+(x_1)^2-(d_2)^2=0,
&&(d_2)^2+(x_2)^2-L^2=0.
\endxalignat
$$
Here $x_1$, $x_2$, $x_3$ are the edges of a cuboid and $d_1$, $d_2$, $d_3$ are 
its face diagonals, while $L$ is its space diagonal. Actually the number of
the equations \mythetag{1.1} can be reduced from six to four since the equations 
of the right column in \mythetag{1.1} are equivalent to one  equation
$(x_1)^2+(x_2)^2+(x_3)^2=L^2$.\par
     Recently in \mycite{40} the equations \mythetag{1.1} were reduced
to a single Diophantine equation with respect to four especially introduced
parameters $a$, $b$, $c$, and $u$. On the base of this equation in \mycite{41}
three cuboid conjectures were formulated. These conjectures were studied in 
\myciterange{42}{42}{--}{44}. However, they are not yet proved.\par
     In \mycite{45} another approach to the equations \mythetag{1.1} was
tested. It is based on the intrinsic $S_3$ symmetry of the equations 
\mythetag{1.1}. Indeed, if $\sigma\in S_3$, then we can write 
$$
\pagebreak 
\xalignat 3
&\hskip -2em
\sigma(x_i)=x_{\sigma i},
&&\sigma(d_i)=d_{\sigma i},
&&\sigma(L)=L.
\mytag{1.2}
\endxalignat
$$
Each transformation $\sigma\in S_3$ permutes the equations \mythetag{1.1}, 
but the system in whole remains unchanged. Factor equations are produced from 
\mythetag{1.1} by introducing new variables which are unchanged under the
transformations \mythetag{1.2}. In \mycite{45} such variables were defined as
values of elementary multisymmetric polynomials:
$$
\gather
\hskip -2em
\aligned
&e_{\sssize [1,0]}=x_1+x_2+x_3,\\
&e_{\sssize [2,0]}=x_1\,x_2+x_2\,x_3+x_3\,x_1,\\
&e_{\sssize [3,0]}=x_1\,x_2\,x_3,
\endaligned
\mytag{1.3}\\
\vspace{2ex}
\hskip -2em
\aligned
&e_{\sssize [0,1]}=d_1+d_2+d_3,\\
&e_{\sssize [0,2]}=d_1\,d_2+d_2\,d_3+d_3\,d_1,\\
&e_{\sssize [0,3]}=d_1\,d_2\,d_3,
\endaligned
\mytag{1.4}\\
\vspace{2ex}
\hskip -2em
\aligned
&e_{\sssize [2,1]}=x_1\,x_2\,d_3+x_2\,x_3\,d_1+x_3\,x_1\,d_2,\\
&e_{\sssize [1,1]}=x_1\,d_2+d_1\,x_2+x_2\,d_3+d_2\,x_3+x_3\,d_1+d_3\,x_1,\\
&e_{\sssize [1,2]}=x_1\,d_2\,d_3+x_2\,d_3\,d_1+x_3\,d_1\,d_2.
\endaligned
\mytag{1.5}
\endgather
$$
The polynomials \mythetag{1.3} coincide with regular elementary symmetric 
polynomials in $x_1$, $x_2$, $x_3$ (see \mycite{46}). The polynomials 
\mythetag{1.4} coincide with regular elementary symmetric polynomials in 
$d_1$, $d_2$, $d_3$. As for the polynomials \mythetag{1.5}, they are actually 
multisymmetric, i\.\,e\. they depend on double set of variables.\par
     General multisymmetric polynomials, which are also known as vector 
symmetric polynomials, diagonally symmetric polynomials, McMahon polynomials 
etc, were initially studied in \myciterange{47}{47}{--}{53} (see also later 
publications \myciterange{54}{54}{--}{67}).\par
     A general multisymmetric polynomial in our case is defined as an element
of the ring $\Bbb Q[x_1,x_2,x_2,d_1,d_2,d_3,L]$ invariant with respect to the
transformations \mythetag{1.2}. The variables $x_1$, $x_2$, $x_3$ and $d_1$, 
$d_2$, $d_3$ are usually arranged into a matrix:
$$
\hskip -2em
M=\Vmatrix x_1 & x_2 &x_3\\
\vspace{1ex}
d_1 & d_2 & d_3\endVmatrix.
\mytag{1.6}
$$
Due to \mythetag{1.2} the group $S_3$ act upon the matrix \mythetag{1.6} by 
permuting its columns. The polynomials from $\Bbb Q[x_1,x_2,x_2,d_1,d_2,d_3,L]$
invariant with respect this action of $S_3$ constitute a ring\footnotemark. 
We denote this ring through $\Sym\!\Bbb Q[M,L]$.\par
\footnotetext{\ Actually both rings $\Bbb Q[x_1,x_2,x_2,d_1,d_2,d_3,L]$ and
$\Sym\!\Bbb Q[M,L]$ are algebras over the field of rational numbers $\Bbb Q$.}
    Let's denote through $p_{\kern 1pt 1}$, $p_{\kern 1pt 2}$, $p_{\kern 1pt 3}$,
$p_{\kern 0.5pt 4}$, $p_{\kern 1pt 5}$, $p_{\kern 1pt 6}$ the left hand sides of 
the cuboid equations \mythetag{1.1}. Then we have the following six polynomials:
$$
\xalignat 2
&\hskip -2em
p_{\kern 1pt 1}=(x_1)^2+(x_2)^2-(d_3)^2,
&&p_{\kern 0.5pt 4}=(d_3)^2+(x_3)^2-L^2,\\
&\hskip -2em
p_{\kern 1pt 2}=(x_2)^2+(x_3)^2-(d_1)^2,
&&p_{\kern 1pt 5}=(d_1)^2+(x_1)^2-L^2,
\mytag{1.7}\\
&\hskip -2em
p_{\kern 1pt 3}=(x_3)^2+(x_1)^2-(d_2)^2,
&&p_{\kern 1pt 6}=(d_2)^2+(x_2)^2-L^2. 
\endxalignat
$$
The polynomials \mythetag{1.7} generate an ideal in the ring 
$\Bbb Q[x_1,x_2,x_2,d_1,d_2,d_3,L]$:
$$
\pagebreak
\hskip -2em
I=\bigl<p_{\kern 1pt 1},p_{\kern 1pt 2},p_{\kern 1pt 3},
p_{\kern 0.5pt 4},p_{\kern 1pt 5},p_{\kern 1pt 6}\bigr>.
\mytag{1.8}
$$
The intersection of the ideal \mythetag{1.8} with the subring $\Sym\!\Bbb Q[M,L]$
of the polynomial ring $\Bbb Q[x_1,x_2,x_2,d_1,d_2,d_3,L]$ is an ideal in
$\Sym\!\Bbb Q[M,L]$: 
$$
\hskip -2em
I_{\text{sym}}=I\cap\Sym\!\Bbb Q[M,L].
\mytag{1.9}
$$
\mydefinition{1.1} A polynomial equation $p\kern 1pt(x_1,x_2,x_2,d_1,d_2,d_3,L)=0$ 
with the polynomial $p\in I_{\text{sym}}$ is called a {\it factor equation\/} of 
the cuboid equations \mythetag{1.1} with respect to their $S_3$ symmetry. 
\enddefinition
     The main goal of this paper is to describe the ideal \mythetag{1.9} in the ring
of multisymmetric polynomials by calculating a finite basis of this ideal.\par
\head
2. The substitution homomorphism. 
\endhead
     Let $\Bbb Q[E,L]=\Bbb Q[E_{10},E_{20},E_{30},E_{01},E_{02},E_{03},E_{21},E_{11},
E_{12},L]$ be a polynomial ring with ten independent variables. If $q\in\Bbb Q[E,L]$,
then substituting the elementary multisymmetric polynomials \mythetag{1.3},
\mythetag{1.4}, and \mythetag{1.5} for $E_{10}$, $E_{20}$, $E_{30}$, $E_{01}$,
$E_{02}$, $E_{03}$, $E_{21}$, $E_{11}$, $E_{12}$ into the arguments of $q$, we get 
a polynomial $p\in\Sym\!\Bbb Q[M,L]$. This means that we have a mapping
$$
\hskip -2em
\varphi\!:\,\Bbb Q[E,L]\longrightarrow\Sym\!\Bbb Q[M,L].
\mytag{2.1}
$$ 
It is easy to see that the mapping \mythetag{2.1} is a ring homomorphism. Such a  
homomorphism is called a {\it substitution homomorphism}. 
\mytheorem{2.1} The elementary multisymmetric polynomials \mythetag{1.3},
\mythetag{1.4}, and \mythetag{1.5} generate the ring of all multisymmetric 
polynomials, i\.\,e\. each multisymmetric polynomial $p\in\Sym\!\Bbb Q[M,L]$ 
can be expressed as a polynomial with rational coefficients through these elementary 
multisymmetric polynomials.
\endproclaim
     The theorem~\mythetheorem{2.1} is known as the fundamental theorem for 
elementary multisymmetric polynomials. Its proof can be found in \mycite{53}. 
The theorem~\mythetheorem{2.1} means that the mapping \mythetag{2.1} is surjective.
Unfortunately the elementary multisymmetric polynomials \mythetag{1.3},
\mythetag{1.4}, and \mythetag{1.5} are not algebraically independent over $\Bbb Q$.
For this reason the homomorphism \mythetag{2.1} is not bijective. It has a nonzero 
kernel:
$$
\hskip -2em
\Ker\varphi=K\neq\{0\}.
\mytag{2.2}
$$
The kernel \mythetag{2.2} is an ideal of the ring $\Bbb Q[E,L]$. According to 
Hilbert's basis theorem (see \mycite{68} and \mycite{69}) each ideal of the 
$\Bbb Q[E,L]$ is finitely generated. This means that 
$$
\hskip -2em
K=\bigl<q_{\kern 1pt 1},\,\ldots,\,q_{\kern 1pt n}\bigr>.
\mytag{2.3}
$$
At present time I know seven polynomials belonging to the ideal \mythetag{2.3}.
They are found by means of direct calculations. Here is the first of these seven 
polynomials:
$$
\pagebreak
\hskip -2em
\gathered
q_{\kern 1pt 1}=
4\,E_{01}\,E_{02}\,E_{20}-E_{02}\,E_{10}^2\,E_{01}-E_{01}^3\,E_{20}\,+\\
+\,E_{10}\,E_{11}\,E_{01}^2-E_{11}^2\,E_{01}-2\,E_{10}\,E_{01}\,E_{12}
+3\,E_{03}\,E_{10}^2\,-\\
-\,9\,E_{03}\,E_{20}-3\,E_{21}\,E_{02}+E_{21}\,E_{01}^2+3\,E_{11}\,E_{12},
\endgathered
\mytag{2.4}
$$
The other six polynomials are given by the following formulas:
$$
\allowdisplaybreaks
\gather
\hskip -2em
\gathered
q_{\kern 1pt 2}=
4\,E_{10}\,E_{20}\,E_{02}-E_{20}\,E_{01}^2\,E_{10}-E_{10}^3\,E_{02}\,+\\
+\,E_{01}\,E_{11}\,E_{10}^2-E_{11}^2\,E_{10}-2\,E_{01}\,E_{10}\,E_{21}
+3\,E_{30}\,E_{01}^2\,-\\
-\,9\,E_{30}\,E_{02}-3\,E_{12}\,E_{20}+E_{12}\,E_{10}^2+3\,E_{11}\,E_{21},
\endgathered
\mytag{2.5}\\
\vspace{2ex}
\hskip -2em
\gathered
q_{\kern 1pt 3}=
9\,E_{21}\,E_{12}-E_{01}^2\,E_{10}\,E_{21}-6\,E_{10}\,E_{11}\,E_{12}
-6\,E_{01}\,E_{12}\,E_{20}\,+\\
+\,5\,E_{01}\,E_{10}^2\,E_{12}-3\,E_{11}^3
+\,7\,E_{10}\,E_{11}^2\,E_{01}+12\,E_{11}\,E_{20}\,E_{02}\,-\\
-\,3\,E_{01}^2\,E_{11}\,E_{20}-3\,E_{02}\,E_{10}^2\,E_{11}
-4\,E_{01}^2\,E_{10}^2\,E_{11}-81\,E_{03}\,E_{30}\,+\\
+\,18\,E_{01}\,E_{02}\,E_{30}-3\,E_{01}^3\,E_{30}+36\,E_{20}\,E_{10}\,E_{03}
-9\,E_{03}\,E_{10}^3\,-\\
-\,16\,E_{01}\,E_{02}\,E_{20}\,E_{10}+4\,E_{01}^3\,E_{10}\,E_{20}
+4\,E_{01}\,E_{10}^3\,E_{02},
\endgathered\quad
\mytag{2.6}\\
\vspace{2ex}
\hskip -2em
\gathered
q_{\kern 1pt 4}=
3\,E_{01}\,E_{21}^2-2\,E_{01}^2\,E_{21}\,E_{20}-9\,E_{01}\,E_{12}\,E_{30}\,+\\
+\,E_{10}\,E_{12}\,E_{01}\,E_{20}
-E_{11}^2\,E_{20}\,E_{01}+3\,E_{01}^2\,E_{30}\,E_{11}
+E_{11}\,E_{20}\,E_{01}^2\,E_{10}\,-\\
-\,3\,E_{01}\,E_{30}\,E_{02}\,E_{10}
+4\,E_{01}\,E_{20}^2\,E_{02}-E_{01}^3\,E_{20}^2-E_{01}\,E_{20}\,E_{10}^2\,E_{02},
\endgathered\quad
\mytag{2.7}\\
\vspace{2ex}
\hskip -2em
\gathered
q_{\kern 1pt 5}=
-27\,E_{10}\,E_{21}\,E_{03}+E_{10}\,E_{01}^3\,E_{21}+9\,E_{10}\,E_{12}^2
-E_{11}^2\,E_{10}\,E_{01}^2\,-\\
-\,6\,E_{02}\,E_{12}\,E_{10}^2-2\,E_{01}^2\,E_{12}\,E_{10}^2
-3\,E_{02}\,E_{11}^2\,E_{10}-E_{01}^2\,E_{10}^3\,E_{02}\,+\\
+\,9\,E_{11}\,E_{03}\,E_{10}^2+3\,E_{01}\,E_{02}\,E_{10}^2\,E_{11}
+E_{01}^3\,E_{11}\,E_{10}^2-3\,E_{10}^3\,E_{02}^2\,+\\
+\,3\,E_{10}^3\,E_{01}\,E_{03}+12\,E_{10}\,E_{20}\,E_{02}^2
+E_{02}\,E_{20}\,E_{01}^2\,E_{10}\,-\\
-\,E_{01}^4\,E_{20}\,E_{10}-18\,E_{10}\,E_{01}\,E_{03}\,E_{20}
+3\,E_{11}\,E_{01}\,E_{10}\,E_{12},
\endgathered\quad
\mytag{2.8}\\
\vspace{2ex}
\hskip -2em
\gathered
q_{\kern 1pt 6}=
-27\,E_{03}\,E_{21}+E_{21}\,E_{01}^3+9\,E_{12}^2
+3\,E_{12}\,E_{01}\,E_{11}\,-\\
-\,2\,E_{01}^2\,E_{10}\,E_{12}-3\,E_{02}\,E_{11}^2-E_{01}^2\,E_{11}^2
+9\,E_{03}\,E_{11}\,E_{10}\,-\\
-\,3\,E_{10}^2\,E_{02}^2+3\,E_{01}\,E_{02}\,E_{11}\,E_{10}
+E_{01}^3\,E_{11}\,E_{10}\,-\\
-18\,E_{20}\,E_{01}\,E_{03}+3\,E_{03}\,E_{01}\,E_{10}^2
-6\,E_{02}\,E_{10}\,E_{12}\,-\\
-\,E_{01}^4\,E_{20}+12\,E_{02}^2\,E_{20}+E_{01}^2\,E_{02}\,E_{20}
-E_{01}^2\,E_{10}^2\,E_{02},
\endgathered\quad
\mytag{2.9}\\
\vspace{2ex}
\hskip -2em
\gathered
q_{\kern 1pt 7}=
3\,E_{21}^2-2\,E_{20}\,E_{01}\,E_{21}-9\,E_{30}\,E_{12}
+E_{10}\,E_{12}\,E_{20}\,-\\
-\,E_{20}\,E_{11}^2+3\,E_{30}\,E_{11}\,E_{01}
+E_{10}\,E_{20}\,E_{11}\,E_{01}\,-\\
-\,3\,E_{02}\,E_{10}\,E_{30}+4\,E_{20}^2\,E_{02}-E_{01}^2\,E_{20}^2
-E_{10}^2\,E_{20}\,E_{02}.
\endgathered\quad
\mytag{2.10}
\endgather
$$
\mytheorem{2.2} Seven polynomials \mythetag{2.4}, \mythetag{2.5},
\mythetag{2.6}, \mythetag{2.7}, \mythetag{2.8}, \mythetag{2.9}, 
\mythetag{2.10} constitute a basis for the ideal $K$ being the kernel 
of the homomorphism \mythetag{2.1}. 
\endproclaim
     Proving the theorem~\mythetheorem{2.1} is an algorithmically 
solvable problem. For this purpose the Gr\"obner bases technique 
should be applied to the ring
$$
\Bbb Q[x_1,x_2,x_3,d_1,d_2,d_3,E_{10},E_{20},E_{30},E_{01},E_{02},E_{03},
E_{21},E_{11},E_{12},L].\quad
\mytag{2.11}
$$
Gr\"obner bases are associated with monomial orderings (see \mycite{69} 
or \mycite{70}). The lexicographic ordering ({\bf lex}) is the most simple 
one. It is defined through some ordering of variables. In the case of the 
ring \mythetag{2.11} one should choose the ordering
$$
\hskip -2em
\gathered
x_1>x_2>x_3>d_1>d_2>d_3>E_{21}>E_{12}>E_{11}>\\
>E_{30}>E_{03}>E_{20}>E_{02}>E_{10}>E_{01}>L.
\endgathered
\mytag{2.12}
$$
Due to the lexicographic ordering based on \mythetag{2.12} each polynomial
$r$ of the ring \mythetag{2.11} gains its leading term $\LT(r)$ with respect
to this {\bf lex}-ordering. 
\mydefinition{2.1} For each ideal $I$ of a polynomial ring the ideal
$\LT(I)$ is generated by leading terms of all polynomials of this ideal.
\enddefinition
\mydefinition{2.2} A basis $r_{\kern 0.5pt 1},\,\ldots,r_s$ of an ideal $I$ 
is called a {\it Gr\"obner basis\/} if the leading terms 
$\LT(r_{\kern 0.5pt 1}),\,\ldots,\LT(r_s)$ generate the ideal $\LT(I)$.
\enddefinition
     An algorithm for computing Gr\"obner bases was first published by 
Bruno Buchberger in 1965 in his PhD thesis \mycite{71}. Wolfgang Gr\"obner 
was Buchberger's thesis adviser. Similar algorithms were developed for local 
rings by Heisuke Hironaka in 1964 (see \mycite{72} and \mycite{73}) and 
for free Lie algebras by A.~I.~Shirshov in 1962 (see \mycite{74}).\par
     The ring \mythetag{2.11} comprises both of the rings $\Bbb Q[M,L]$ and
$\Bbb Q[E,L]$. For this reason one can consider the following nine polynomials
in this ring:
$$
\xalignat 3
&\hskip -2em
r_{\kern 0.5pt 1}=E_{10}-e_{\sssize [1,0]},
&&r_{\kern 0.5pt 2}=E_{20}-e_{\sssize [2,0]},
&&r_{\kern 0.5pt 3}=E_{30}-e_{\sssize [3,0]},
\qquad\\
&\hskip -2em
r_{\kern 0.5pt 4}=E_{01}-e_{\sssize [0,1]},
&&r_{\kern 0.5pt 5}=E_{02}-e_{\sssize [0,2]},
&&r_{\kern 0.5pt 6}=E_{03}-e_{\sssize [0,3]},
\qquad
\mytag{2.13}\\
&\hskip -2em
r_{\kern 0.5pt 7}=E_{21}-e_{\sssize [2,1]},
&&r_{\kern 0.5pt 8}=E_{11}-e_{\sssize [1,1]},
&&r_{\kern 0.5pt 9}=E_{12}-e_{\sssize [1,2]}.
\qquad
\endxalignat
$$
The polynomials \mythetag{2.13} are constructed with the use of the
elementary multisymmetric polynomials \mythetag{1.3}, \mythetag{1.4}, 
and \mythetag{1.5}. They generate the ideal 
$$
\hskip -2em
K_0=\bigl<r_{\kern 0.5pt 1},r_{\kern 0.5pt 2},r_{\kern 0.5pt 3},
r_{\kern 0.5pt 4},r_{\kern 0.5pt 5},r_{\kern 0.5pt 6},
r_{\kern 0.5pt 7},r_{\kern 0.5pt 8},r_{\kern 0.5pt 9}\bigr>
\mytag{2.14}
$$
of the ring \mythetag{2.11}. The kernel of the homomorphism \mythetag{2.1} 
in \mythetag{2.2} coincides with the $6$-th elimination ideal for the ideal 
\mythetag{2.14} with respect to the ordering \mythetag{2.12}:
$$
K=\Ker\varphi=K_6=K_0\cap\Bbb Q[E,L].
\mytag{2.15}
$$ 
\mydefinition{2.3} Let $I$ be an ideal in the polynomial ring $\Bbb Q[x_1,
\ldots,x_n]$. Then the intersection of the ideal $I$ with the subring 
$\Bbb Q[x_{k+1},\ldots,x_n]\subset\Bbb Q[x_1,\ldots,x_n]$ is called the
{\it $k$-th elimination ideal\/} of the ideal $I$:
$$
\hskip -2em
I_k=I\cap \Bbb Q[x_{k+1},\ldots,x_n].
\mytag{2.16}
$$ 
\enddefinition
\mytheoremwithtitle{2.3}{ (elimination theorem)} Let $I$ be an ideal in the 
polynomial ring $\Bbb Q[x_1,\ldots,x_n]$ and let $G=\{g_1,\,\ldots,\,g_s\}$ be
its Gr\"obner basis with respect to the {\bf lex}-ordering with
$x_1>x_2>\ \ldots\ >x_n$. Then for any $0\leqslant k\leqslant n$ the intersection
$$
\hskip -2em
G_k=G\cap \Bbb Q[x_{k+1},\ldots,x_n]
\mytag{2.17}
$$
is a Gr\"obner basis for the $k$-th elimination ideal $I_k$.
\endproclaim
The definition~\mythedefinition{2.3} and the formula \mythetag{2.16} explain
the formula \mythetag{2.15}, while the theorem~\mythetheorem{2.3} along with
the formula \mythetag{2.17} yields an algorithm for calculating a basis for
the ideal \mythetag{2.2} and thus for proving the theorem~\mythetheorem{2.2}. 
The proof of the theorem~\mythetheorem{2.3} can be found in \mycite{69}.\par
     The algorithm provided by the theorem~\mythetheorem{2.3} is already 
implemented in many packages for symbolic computations. For instance, the
Maxima package (version 5.22.1) contains the Gr\"obner subpackage 
(revision 1.6) with the command 
$$
\hskip -2em
\text{\tt poly\_\,elimination\,\_\,ideal(L,k,V)},
\mytag{2.18}
$$
where $L$ is a list of polynomials, $k$ is the integer number from \mythetag{2.17},
and $V$ is a list of variables. Due to \mythetag{2.12}, \mythetag{2.13}, and
\mythetag{2.15} in my case I have $k=6$ and
$$
\aligned
&L=[r_{\kern 0.5pt 1},r_{\kern 0.5pt 2},r_{\kern 0.5pt 3},r_{\kern 0.5pt 4},
r_{\kern 0.5pt 5},r_{\kern 0.5pt 6},r_{\kern 0.5pt 7},r_{\kern 0.5pt 8},
r_{\kern 0.5pt 9}],\\
&V=[x_1,x_2,x_3,d_1,d_2,d_3,E_{21},E_{12},E_{11},
E_{30},E_{03},E_{20},E_{02},E_{10},E_{01},L].
\endaligned
$$ 
After running the command \mythetag{2.18} with the above parameters on a machine 
with dual core Prescott 2.8E Intel Pentium-4 processor and with 500 megabytes RAM 
on board I have got a Gr\"obner basis $G_K$ of the ideal $K$ consisting of $14$ 
polynomials. Some of them are rather huge for to typeset them here. Using this 
Gr\"obner basis, I have verified that the polynomials \mythetag{2.4}, \mythetag{2.5},
\mythetag{2.6}, \mythetag{2.7}, \mythetag{2.8}, \mythetag{2.9}, and \mythetag{2.10}
do actually belong to the kernel of the homomorphism \mythetag{2.1}.\par
     Conversely, the polynomials \mythetag{2.4}, \mythetag{2.5}, \mythetag{2.6}, 
\mythetag{2.7}, \mythetag{2.8}, \mythetag{2.9}, and \mythetag{2.10} generate 
their own Gr\"obner basis $G_Q$. Using this second Gr\"obner basis $G_Q$, I have
tested each polynomial of the first Gr\"obner basis $G_K$ and have found that all
of these polynomials belong to the ideal $Q=\bigl<q_{\kern 1pt 1},q_{\kern 1pt 2},
q_{\kern 1pt 3},q_{\kern 1pt 4},q_{\kern 1pt 5},q_{\kern 1pt 6},q_{\kern 1pt 7}
\bigr>$ generated by the polynomials \mythetag{2.4}, \mythetag{2.5}, \mythetag{2.6}, 
\mythetag{2.7}, \mythetag{2.8}, \mythetag{2.9}, and \mythetag{2.10}. This result
means that the ideals $K=\Ker\varphi$ and $Q$ do coincide, i\.\,e\. I have got a
computer aided proof of the theorem~\mythetheorem{2.2}.\par
\head
3. The fine structure of the ideal $I_{\text{sym}}$. 
\endhead
     The ideal $I$ producing $I_{\text{sym}}$ in \mythetag{1.9}
is generated by six polynomials \mythetag{1.7} in \mythetag{1.8}. Actually,
the number of generating polynomials of the ideal $I$ can be reduced from 
six to four. Indeed, we can write
$$
\hskip -2em
I=\bigl<p_{\kern 1pt 0},p_{\kern 1pt 1},p_{\kern 1pt 2},p_{\kern 1pt 3}\bigr>,
\mytag{3.1}
$$
where $p_{\kern 1pt 0}$ is a symmetric polynomial given by the formula
$$
\hskip -2em
p_{\kern 1pt 0}=(x_1)^2+(x_2)^2+(x_3)^2-L^2.
\mytag{3.2}
$$
Due to the relationship \mythetag{3.1} each polynomial $p\in I_{\text{sym}}$ 
is written as 
$$
\hskip -2em
p=\alpha_0\,p_{\kern 1pt 0}+\sum^3_{i=1}\alpha_i\,p_{\kern 1pt i},
\mytag{3.3}
$$
where $\alpha_i\in\Bbb Q[M,L]$. Since $p$ is a multisymmetric polynomial, it
should be invariant with respect to the symmetrization operator $S$ defined
by the formula
$$
\hskip -2em
S(p)=\sum_{\sigma\in S_3}\frac{\sigma^{-1}(p)}{6}.
\mytag{3.4}
$$
The invariance of $p$ with respect to the operator \mythetag{3.4} is written
as $p=S(p)$. Therefore, applying $S$ to \mythetag{3.3}, we derive the formula
$$
\hskip -2em
p=S(\alpha_0\,p_{\kern 1pt 0})+\sum^3_{i=1}S(\alpha_i\,p_{\kern 1pt i}).
\mytag{3.5}
$$\par
     Now let's recall the formulas \mythetag{1.2}. Applying them to the 
polynomials \mythetag{1.7} and \mythetag{3.2}, we derive the analogous 
formulas 
$$
\xalignat 2
&\hskip -2em
\sigma(p_{\kern 1pt i})=p_{\kern 1pt\sigma i},
&&\sigma(p_{\kern 1pt 0})=p_{\kern 1pt 0}
\mytag{3.6}
\endxalignat
$$
for $p_{\kern 1pt 1}$, $p_{\kern 1pt 2}$, $p_{\kern 1pt 3}$, 
and $p_{\kern 1pt 0}$. Relying on \mythetag{3.6} we introduce the following 
notations:
$$
\xalignat 2
\hskip -2em
&\tilde\alpha_0=S(\alpha_0),
&&\tilde\alpha_i=\sum_{\sigma\in S_3}\frac{\sigma^{-1}(\alpha_{\sigma i})}{6}
\mytag{3.7}
\endxalignat
$$
Using \mythetag{3.7}, we can transform the formula \mythetag{3.5} as follows:
$$
\hskip -2em
p=\tilde\alpha_0\,p_{\kern 1pt 0}+\sum^3_{i=1}\tilde\alpha_i\,p_{\kern 1pt i}.
\mytag{3.8}
$$
The formula \mythetag{3.8} is analogous to the formula \mythetag{3.3}. However,
unlike the original coefficients $\alpha_1$, $\alpha_2$, $\alpha_3$, and
$\alpha_0$ in \mythetag{3.3}, the coefficients \mythetag{3.7} obey the
relationships
$$
\xalignat 2
&\hskip -2em
\sigma(\tilde\alpha_i)=\tilde\alpha_{\sigma i},
&&\sigma(\tilde\alpha_0)=\tilde\alpha_0.
\mytag{3.9}
\endxalignat
$$
The formulas \mythetag{3.8} and \mythetag{3.9} mean that we have proved
the following lemma.
\mylemma{3.1} Each polynomial $p\in I_{\text{sym}}=I\cap\Sym\!\Bbb Q[M,L]$ 
is given by the formula \mythetag{3.3} with the coefficients
$\alpha_i\in\Bbb Q[M,L]$ obeying the relationships
$$
\xalignat 2
&\hskip -2em
\sigma(\alpha_i)=\alpha_{\sigma i},
&&\sigma(\alpha_0)=\alpha_0.
\mytag{3.10}
\endxalignat
$$
\endproclaim
     The formulas \mythetag{3.10} in the lemma~\mythelemma{3.1} are
important since, applying them back to the formula \mythetag{3.5} and
taking into account \mythetag{3.6}, we derive
$$
\hskip -2em
p=\alpha_0\,p_{\kern 1pt 0}+3\,S(\alpha_1\,p_{\kern 1pt 1}).
\mytag{3.11}
$$
Note that $\alpha_1\in\Bbb Q[M,L]$ in \mythetag{3.11} is a polynomial,
i\.\,e\. it is a sum of monomials:
$$
\hskip -2em
\alpha_1=\!\!\!\!\sum\Sb i,j,k\\m,n,r,s\endSb\!\!\!
\theta_{ijkmnrs}\ x_1^i\,x_2^j\,x_3^k\,d_1^{\kern 1pt m}
\,d_2^{\kern 1pt n}\,d_3^{\kern 1pt r}\,L^s.
\mytag{3.12}
$$
Substituting \mythetag{3.12} into the formula \mythetag{3.11}, we easily
derive the following lemma.
\mylemma{3.2} The ideal $I_{\text{sym}}=I\cap\Sym\!\Bbb Q[M,L]$ of the
ring $\Sym\!\Bbb Q[M,L]$ is gene\-rated by the polynomial $p_{\kern 1pt 0}$ 
and by various polynomials of the form
$$
\hskip -2em
S(p_{\kern 1pt 1}\,x_1^i\,x_2^j\,x_3^k\,d_1^{\kern 1pt m}
\,d_2^{\kern 1pt n}\,d_3^{\kern 1pt r}\,L^s).
\mytag{3.13}
$$
\endproclaim
     Note that the factor $L^s$ is invariant with respect to the operator
$S$. It can be split out from the polynomial \mythetag{3.13}. Similarly, if 
$\mu=\min(i,j,k)>0$ and/or $\nu=\min(m,n,r)>0$, we can split out the invariant
factors $(x_1\,x_2\,x_3)^\mu$  and/or $(d_1\,d_2\,d_3)^\nu$. \pagebreak As a 
result we modify the lemma~\mythelemma{3.2} as follows.
\mylemma{3.3} The ideal $I_{\text{sym}}=I\cap\Sym\!\Bbb Q[M,L]$ of the
ring $\Sym\!\Bbb Q[M,L]$ is gene\-rated by the polynomial $p_{\kern 1pt 0}$ 
and by various polynomials of the form
$$
S(p_{\kern 1pt 1}\,x_1^i\,x_2^j\,x_3^k\,d_1^{\kern 1pt m}
\,d_2^{\kern 1pt n}\,d_3^{\kern 1pt r}),
$$
where at least one of the nonnegative numbers $i,\,j,\,k$ is zero and at least 
one of the nonnegative numbers $m,\,n,\,r$ is zero.
\endproclaim
     The lemma~\mythelemma{3.3} yields a basis for the ideal $I_{\text{sym}}$.
However this basis is not finite. Getting a finite basis of the ideal 
$I_{\text{sym}}$ is a little bit more tricky.\par
\head
4. Partially multisymmetric polynomials.
\endhead 
     Let's consider the formulas \mythetag{3.10}. The polynomial $\alpha_0$
in \mythetag{3.10} is multisymmetric, i\.\,e\. it is invariant with respect 
to the transformations \mythetag{1.2} for all $\sigma\in S_3$. As for the 
polynomials $\alpha_1$, $\alpha_2$, and $\alpha_3$ in \mythetag{3.10}, they
are partially multisymmetric. The formulas \mythetag{3.10} for these 
polynomials yield 
$$
\align
&\hskip -2em
\sigma(\alpha_1)=\alpha_1\text{\ \ if and only if \ }\sigma 1=1,\\
&\hskip -2em
\sigma(\alpha_2)=\alpha_2\text{\ \ if and only if \ }\sigma 2=2,
\mytag{4.1}\\
&\hskip -2em
\sigma(\alpha_3)=\alpha_3\text{\ \ if and only if \ }\sigma 3=3.
\endalign
$$
The formulas \mythetag{4.1} mean that the polynomials $\alpha_1$, $\alpha_2$, 
and $\alpha_3$ are $S_2$ invariant, but they are invariant with respect to
different subgroups of the group $S_3$ isomorphic to the group $S_2$. In order
to describe such partially multisymmetric polynomials we split out the following
three matrices from the matrix \mythetag{1.6}:
$$
\xalignat 3
&M_1=\Vmatrix x_2 &x_3\\ \vspace{1ex} d_2 & d_3\endVmatrix,
&&M_2=\Vmatrix x_1 &x_3\\ \vspace{1ex} d_1 & d_3\endVmatrix,
&&M_3=\Vmatrix x_1 &x_2\\ \vspace{1ex} d_1 & d_2\endVmatrix.
\qquad\quad
\mytag{4.2}
\endxalignat
$$
Like the matrix \mythetag{1.6}, the matrices \mythetag{4.2} can be used for
producing elementary multisymmetric polynomials. Here are these polynomials:
$$
\xalignat 2
&\hskip -2em
f_{\sssize [1,0]}[1]=x_2+x_3,
&&f_{\sssize [2,0]}[1]=x_2\,x_3,
\quad\\
&\hskip -2em
f_{\sssize [0,1]}[1]=d_2+d_3,
&&f_{\sssize [0,1]}[1]=d_2\,d_3,
\quad\\
&\hskip -2em
f_{\sssize [1,0]}[2]=x_3+x_1,
&&f_{\sssize [2,0]}[2]=x_3\,x_1,
\quad\\
\vspace{-1.7ex}
\mytag{4.3}\\
\vspace{-1.7ex}
&\hskip -2em
f_{\sssize [0,1]}[2]=d_3+d_1,
&&f_{\sssize [0,1]}[2]=d_3\,d_1,
\quad\\
&\hskip -2em
f_{\sssize [1,0]}[3]=x_1+x_2,
&&f_{\sssize [2,0]}[3]=x_1\,x_2,
\quad\\
&\hskip -2em
f_{\sssize [0,1]}[3]=d_1+d_2,
&&f_{\sssize [0,1]}[3]=d_1\,d_2.
\quad
\endxalignat
$$
Apart from \mythetag{4.3} there are three other elementary
multisymmetric polynomials:
$$
\align
&\hskip -2em
f_{\sssize [1,1]}[1]=x_2\,d_3+x_3\,d_2,\\
&\hskip -2em
f_{\sssize [1,1]}[2]=x_3\,d_1+x_1\,d_3,
\mytag{4.4}\\
&\hskip -2em
f_{\sssize [1,1]}[3]=x_1\,d_2+x_2\,d_1.
\endalign
$$
The polynomials in \mythetag{4.3} and \mythetag{4.4} are subdivided into 
three groups depending on which matrix \mythetag{4.2} is used for
their production.\par
     Like $\alpha_1$, $\alpha_2$, and $\alpha_3$, the polynomials 
\mythetag{4.3} and \mythetag{4.4} are partially multisymmetric.
They obey the following relationships very similar to \mythetag{3.10}:
$$
\xalignat 2
&\hskip -2em
\sigma(f_{\sssize [1,0]}[i])=f_{\sssize [1,0]}[\sigma i],
&&\sigma(f_{\sssize [2,0]}[i])=f_{\sssize [2,0]}[\sigma i],\\
&\hskip -2em
\sigma(f_{\sssize [0,1]}[i])=f_{\sssize [0,1]}[\sigma i],
&&\sigma(f_{\sssize [0,2]}[i])=f_{\sssize [0,2]}[\sigma i],
\mytag{4.5}\\
&\hskip -2em
\sigma(f_{\sssize [1,1]}[i])=f_{\sssize [1,1]}[\sigma i].
\endxalignat
$$
The polynomials \mythetag{4.3} and \mythetag{4.4} obey a theorem similar
to the theorem~\mythetheorem{2.1}. 
\mytheorem{4.1} The elementary multisymmetric polynomials $f_{\sssize [1,0]}$,
$f_{\sssize [2,0]}$, $f_{\sssize [0,1]}$, $f_{\sssize [0,2]}$, $f_{\sssize [1,1]}$
generate the ring of all $S_2$ multisymmetric polynomials, i\.\,e\. each 
$S_2$ multisymmetric polynomial can be expressed as a polynomial with rational 
coefficients through these elementary multisymmetric polynomials.
\endproclaim
     The theorem~\mythetheorem{4.1} is an $S_2$ version of the fundamental
theorem on elementary multisymmetric polynomials which is formulated for the
general case of $S_n$ multisymmetric polynomials (see \mycite{53}). Applying 
this theorem to $\alpha_1$, $\alpha_2$, and $\alpha_3$, we \nolinebreak get
$$
\hskip -2em
\alpha_i=q_i(x_i,d_i,f_{\sssize [1,0]}[i],f_{\sssize [2,0]}[i],
f_{\sssize [0,1]}[i],f_{\sssize [0,2]}[i],f_{\sssize [1,1]}[i],L),
\mytag{4.6}
$$
where $q_i$ is some polynomial of eight independent variables. The polynomials
$\alpha_1$, $\alpha_2$, and $\alpha_3$ are not independent. They are related
to each other by means of the formulas \mythetag{3.10}. Therefore, applying
\mythetag{1.2}, \mythetag{3.10}, and \mythetag{4.5} to \mythetag{4.6}, we 
conclude that the polynomials $q_i$ in \mythetag{4.6} can be chosen so that
they do coincide, i\.\,e\.
$$
q_1=q_2=q_3=q(x,d,f_{\sssize [1,0]},f_{\sssize [2,0]},
f_{\sssize [0,1]},f_{\sssize [0,2]},f_{\sssize [1,1]},L)
\mytag{4.7}
$$
Applying \mythetag{4.7} to \mythetag{4.6}, we write \mythetag{4.6} as follows:
$$
\hskip -2em
\alpha_i=q(x_i,d_i,f_{\sssize [1,0]}[i],f_{\sssize [2,0]}[i],
f_{\sssize [0,1]}[i],f_{\sssize [0,2]}[i],f_{\sssize [1,1]}[i],L),
\mytag{4.8}
$$\par
     Now let's return back to the formulas \mythetag{3.5} and \mythetag{3.11}. 
Applying \mythetag{4.8} to \mythetag{3.11}, we get the following
expression for $p$:
$$
\hskip -2em
p=\alpha_0\,p_{\kern 1pt 0}+3\,S(q\,p_{\kern 1pt 1}).
\mytag{4.9}
$$
Here $q=q(x_1,d_1,f_{\sssize [1,0]}[1],f_{\sssize [2,0]}[1],
f_{\sssize [0,1]}[1],f_{\sssize [0,2]}[1],f_{\sssize [1,1]}[1],L)$ and 
$S$ is the symmetri\-zation operator \mythetag{3.4}. The formula \mythetag{4.9} 
applies to any polynomial $p\in I_{\text{sym}}$.\par
\head
5. The module structure of the ideal $I_{\text{sym}}$.
\endhead
     Each ideal is a module over that ring for which it is an ideal. When applied 
to the ideal $I_{\text{sym}}$, this fact means that
$$
\hskip -2em
p\in I_{\text{sym}}\text{\ \ implies \ }\alpha\,p\in I_{\text{sym}}
\text{\ \ for any \ }\alpha\in\Sym\!\Bbb Q[M,L].
\mytag{5.1} 
$$
Relying on \mythetag{5.1}, let us consider the the product $\alpha\,p$ for a
polynomial $p$ given by the formula \mythetag{4.9}. As a result we obtain the
formula
$$
\alpha p=\alpha\,\alpha_0\,p_{\kern 1pt 0}+3\,\alpha\,S(q\,p_{\kern 1pt 1}).
\quad
\mytag{5.2}
$$
Note that $\alpha$ in \mythetag{5.2} is a multisymmetric polynomial. Therefore 
it goes through the symmetrization operator $S$ as a scalar factor. This yields
$$
\alpha p=\alpha\,\alpha_0\,p_{\kern 1pt 0}+3\,S(\alpha\,q\,p_{\kern 1pt 1}).
\quad
\mytag{5.3}
$$
Comparing the formulas \mythetag{5.3} and \mythetag{4.9}, we conclude that the 
multiplication by $\alpha$ in $I_{\text{sym}}$ is equivalent to the transformation
$$
\xalignat 2
&\hskip -2em
\alpha_0\mapsto\alpha\,\alpha_0,
&&q\mapsto\alpha\,q.
\mytag{5.4}
\endxalignat
$$
The polynomial $\alpha$ in the formulas \mythetag{5.4} is expressed the through 
elementary multisymmetric polynomials \mythetag{1.3}, \mythetag{1.4}, and 
\mythetag{1.5}:
$$
\hskip -2em
\alpha=\alpha(e_{\sssize [1,0]},e_{\sssize [2,0]},e_{\sssize [3,0]},
e_{\sssize [0,1]},e_{\sssize [0,2]},e_{\sssize [0,3]},e_{\sssize [2,1]},
e_{\sssize [1,1]},e_{\sssize [1,2]},L),
\mytag{5.5}
$$
while the polynomial $q$ in \mythetag{5.4} is given by the formula 
$$
\hskip -2em
q=q(x_1,d_1,f_{\sssize [1,0]}[1],f_{\sssize [2,0]}[1],
f_{\sssize [0,1]}[1],f_{\sssize [0,2]}[1],f_{\sssize [1,1]}[1],L)
\mytag{5.6}
$$
\par
     Formally, the polynomials $\alpha$ and $q$ depend on different sets of 
variables, though due to \mythetag{1.3}, \mythetag{1.4}, \mythetag{1.5},
\mythetag{4.3}, and \mythetag{4.4} both sets reduce to $x_1$, $x_2$, $x_3$, 
$d_1$, $d_2$, $d_3$, and $L$. Our next goal is to study the mutual relations 
of arguments in \mythetag{5.5} and \mythetag{5.6}. By means of direct 
calculations we derive the formulas
$$
\xalignat 2
&\hskip -2em
f_{\sssize [1,0]}[1]=e_{\sssize [1,0]}-x_1,
&&f_{\sssize [2,0]}[1]=e_{\sssize [2,0]}-x_1\,e_{\sssize [1,0]}+x_1^2,
\quad\\
\vspace{-1.7ex}
\mytag{5.7}\\
\vspace{-1.7ex}
&\hskip -2em
f_{\sssize [0,1]}[1]=e_{\sssize [0,1]}-x_1,
&&f_{\sssize [0,2]}[1]=e_{\sssize [0,2]}-d_1\,e_{\sssize [0,1]}
+d_1^{\kern 1pt 2}.
\quad
\endxalignat
$$
The polynomial $f_{\sssize [1,1]}[1]$ is reexpressed by the formula 
$$
\hskip -2em
f_{\sssize [1,1]}[1]=e_{\sssize [1,1]}-d_1\,e_{\sssize [1,0]}
-x_1\,e_{\sssize [0,1]}+2\,d_1\,x_1.
\mytag{5.8}
$$
In addition to \mythetag{5.7} and \mythetag{5.8}, there are the following
four equations:
$$
\gather
\hskip -2em
\aligned
x_1^3=x_1^2\,e_{\sssize [1,0]}-x_1\,e_{\sssize [2,0]}+e_{\sssize [3,0]},\\
d_1^{\kern 1pt 3}=d_1^{\kern 1pt 2}\,e_{\sssize [0,1]}-d_1\,e_{\sssize [0,2]}
+e_{\sssize [0,3]},
\endaligned
\mytag{5.9}\\
\vspace{2ex}
\hskip -2em
\aligned
&d_1\,x_1^2=\frac{2\,d_1\,x_1}{3}\,e_{\sssize [1,0]}
+\frac{x_1^2}{3}\,e_{\sssize [0,1]}
-\frac{x_1}{3}\,e_{\sssize [1,1]}
-\frac{d_1}{3}\,e_{\sssize [2,0]}
+\frac{1}{3}\,e_{\sssize [2,1]},\\
&x_1\,d_1^{\kern 1pt 2}=\frac{2\,x_1\,d_1}{3}\,e_{\sssize [0,1]}
+\frac{d_1^{\kern 1pt 2}}{3}\,e_{\sssize [1,0]}
-\frac{d_1}{3}\,e_{\sssize [1,1]}
-\frac{x_1}{3}\,e_{\sssize [0,2]}
+\frac{1}{3}\,e_{\sssize [1,2]}.
\endaligned
\mytag{5.10}
\endgather
$$
The equations \mythetag{5.9} and \mythetag{5.10} are easily derived by
means of direct calculations with the use of the formulas \mythetag{1.3}, 
\mythetag{1.4}, and \mythetag{1.5}.\par
     Let's substitute \mythetag{5.7} and \mythetag{5.8} into the arguments
of the polynomial \mythetag{5.6}. As a result the polynomial $q$ is expressed 
in the form
$$
\hskip -2em
q=\tilde q(x_1,d_1,e_{\sssize [1,0]},e_{\sssize [2,0]},e_{\sssize [3,0]},
e_{\sssize [0,1]},e_{\sssize [0,2]},e_{\sssize [0,3]},e_{\sssize [2,1]},
e_{\sssize [1,1]},e_{\sssize [1,2]},L),
\mytag{5.11}
$$
where $\tilde q$ is some arbitrary polynomial of twelve variables. The first 
formula \mythetag{5.9} expresses $x_1^3$ through $x_1^2$ and $x_1$. Similarly,
the second formula \mythetag{5.9} expresses $d_1^{\kern 1pt 3}$ through 
$d_1^{\kern 1pt 2}$ and $d_1$. Therefore, without loss of generality we can
assume that the order of the polynomial $\tilde q$ in $x_1$ and in $d_1$ is
not higher than $2$, i\.\,e\. the variables $x_1$ and $d_1$ enter this polynomial
through the following monomials:
$$
\xalignat 9
&x_1^2\,d_1^{\kern 1pt 2},
&&x_1^2\,d_1,
&&x_1\,d_1^{\kern 1pt 2},
&&x_1\,d_1,
&&x_1^2,
&&d_1^{\kern 1pt 2},
&&x_1,
&&d_1,
&&1.
\quad\qquad
\mytag{5.12}
\endxalignat
$$
Due to the equations \mythetag{5.10} we can exclude the monomials $x_1^2\,d_1$
and $x_1\,d_1^{\kern 1pt 2}$ from the above list \mythetag{5.12} and write
the formula \mythetag{5.11} as
$$
q=Q_{22}\,x_1^2\,d_1^{\kern 1pt 2}+Q_{11}\,x_1\,d_1+Q_{20}\,x_1^2
+Q_{02}\,d_1^{\kern 1pt 2}+Q_{10}\,x_1+Q_{01}\,d_1+Q_{00}.
\quad
\mytag{5.13}
$$
The coefficients $Q_{ij}$ in \mythetag{5.13} are produced by polynomials of ten 
variables:
$$
\hskip -2em
Q_{ij}=Q_{ij}(e_{\sssize [1,0]},e_{\sssize [2,0]},e_{\sssize [3,0]},
e_{\sssize [0,1]},e_{\sssize [0,2]},e_{\sssize [0,3]},e_{\sssize [2,1]},
e_{\sssize [1,1]},e_{\sssize [1,2]},L).
\mytag{5.14}
$$
\par
     The values of the expressions \mythetag{5.5} and \mythetag{5.14} are regular
multisymmetric polynomials from the ring $\Sym\!\Bbb Q[M,L]$, while the values of
the expression \mythetag{5.11} constitute a module over this ring. Due to 
\mythetag{5.13} this module is finitely generated. 
\head
6. A basis of the ideal $I_{\text{sym}}$.
\endhead
     Now we can substitute the formula \mythetag{5.13} with the coefficients
\mythetag{5.14} into the formula \mythetag{4.9}. As a result we can write 
\mythetag{4.9} as 
$$
\gathered
p=\alpha_0\,p_{\kern 1pt 0}
+Q_{22}\,S(3\,x_1^2\,d_1^{\kern 1pt 2}\,p_{\kern 1pt 1})
+Q_{11}\,S(3\,x_1\,d_1\,p_{\kern 1pt 1})
+Q_{20}\,S(3\,x_1^2\,p_{\kern 1pt 1})\,+\\
+\,Q_{02}\,S(3\,d_1^{\kern 1pt 2}\,p_{\kern 1pt 1})
+Q_{10}\,S(3\,x_1\,p_{\kern 1pt 1})
+Q_{01}\,S(3\,d_1\,p_{\kern 1pt 1})
+Q_{00}\,S(3\,p_{\kern 1pt 1}),
\endgathered\qquad
\mytag{6.1}
$$
where $p$ is an arbitrary polynomial from the ideal $I_{\text{sym}}$ and
$\alpha_0$, $Q_{22}$, $Q_{11}$, $Q_{10}$, $Q_{01}$, $Q_{00}$ are arbitrary 
polynomials from the ring $\Sym\!\Bbb Q[M,L]$. The formula \mythetag{6.1}
proves the following theorem, which is the main result of the present paper.
\mytheorem{6.1} The ideal $I_{\text{sym}}$ in the ring $\Sym\!\Bbb Q[M,L]$
defined by the left hand sides of the cuboid equations \mythetag{1.1}
through the formulas \mythetag{1.7}, \mythetag{1.8}, \mythetag{1.9} is 
finitely generated. Eight multisymmetric polynomials
$$
\xalignat 2
&\hskip -2em
\tilde p_{\kern 1pt 1}=p_{\kern 1pt 0},
&&\tilde p_{\kern 1pt 2}=S(3\,p_{\kern 1pt 1}),\\
&\hskip -2em
\tilde p_{\kern 1pt 3}=S(3\,d_1\,p_{\kern 1pt 1}),
&&\tilde p_{\kern 1pt 4}=S(3\,x_1\,p_{\kern 1pt 1}),\\
\vspace{-1.7ex}
\mytag{6.2}\\
\vspace{-1.7ex}
&\hskip -2em
\tilde p_{\kern 1pt 5}=S(3\,x_1\,d_1\,p_{\kern 1pt 1}),
&&\tilde p_{\kern 1pt 6}=S(3\,x_1^2\,p_{\kern 1pt 1}),\\
&\hskip -2em
\tilde p_{\kern 1pt 7}=S(3\,d_1^{\kern 1pt 2}\,p_{\kern 1pt 1}),
&&\tilde p_{\kern 1pt 8}=S(3\,x_1^2\,d_1^{\kern 1pt 2}\,p_{\kern 1pt 1})
\endxalignat
$$ 
belong to the ideal $I_{\text{sym}}$ and constitute a basis of this ideal.
\endproclaim
     The polynomial $\tilde p_{\kern 1pt 1}=p_{\kern 1pt 0}$ from the first 
formula \mythetag{6.2} is already known in an explicit form. It is given by 
the formula \mythetag{3.2}. The polynomial $p_1$ used in the other formulas 
\mythetag{6.2} is also known in an explicit form (see \mythetag{1.7}). Now, 
applying the formula \mythetag{3.4} for $S$, we can explicitly calculate the 
polynomials $\tilde p_{\kern 1pt 2}$, $\tilde p_{\kern 1pt 3}$,
$\tilde p_{\kern 1pt 4}$, $\tilde p_{\kern 1pt 5}$, $\tilde p_{\kern 1pt 6}$,
$\tilde p_{\kern 1pt 7}$, and $\tilde p_{\kern 1pt 8}$. Here is the formula
for the polynomial $\tilde p_{\kern 1pt 2}$:
$$
\hskip -2em
\tilde p_{\kern 1pt 2}=(x_2^2+x_3^2-d_1^{\kern 1pt 2})
+(x_3^2+x_1^2-d_2^{\kern 1pt 2})+(x_1^2+x_2^2-d_3^{\kern 1pt 2}).
\mytag{6.3}
$$
The explicit formulas for $\tilde p_{\kern 1pt 3}$,
$\tilde p_{\kern 1pt 4}$, $\tilde p_{\kern 1pt 5}$, $\tilde p_{\kern 1pt 6}$,
$\tilde p_{\kern 1pt 7}$, and $\tilde p_{\kern 1pt 8}$ are listed just below:
$$
\gather
\hskip -2em
\tilde p_{\kern 1pt 3}=d_1\,(x_2^2+x_3^2-d_1^{\kern 1pt 2})
+d_2\,(x_3^2+x_1^2-d_2^{\kern 1pt 2})+d_3\,(x_1^2+x_2^2-d_3^{\kern 1pt 2}),
\quad
\mytag{6.4}\\
\hskip -2em
\tilde p_{\kern 1pt 4}=x_1\,(x_2^2+x_3^2-d_1^{\kern 1pt 2})
+x_2\,(x_3^2+x_1^2-d_2^{\kern 1pt 2})+x_3\,(x_1^2+x_2^2-d_3^{\kern 1pt 2}),
\quad
\mytag{6.5}\\
\hskip -9em
\gathered
\tilde p_{\kern 1pt 5}=x_1\,d_1\,(x_2^2+x_3^2-d_1^{\kern 1pt 2})
+x_2\,d_2\,(x_3^2+x_1^2-d_2^{\kern 1pt 2})\,+\\
\vphantom{1}\hskip 23em +\,x_3\,d_3\,(x_1^2+x_2^2-d_3^{\kern 1pt 2}),
\endgathered
\mytag{6.6}\\
\hskip -2em
\tilde p_{\kern 1pt 6}=x_1^2\,(x_2^2+x_3^2-d_1^{\kern 1pt 2})
+x_2^2\,(x_3^2+x_1^2-d_2^{\kern 1pt 2})+x_3^2\,(x_1^2+x_2^2-d_3^{\kern 1pt 2}),
\quad
\mytag{6.7}\\
\hskip -2em
\tilde p_{\kern 1pt 7}=d_1^{\kern 1pt 2}\,(x_2^2+x_3^2-d_1^{\kern 1pt 2})
+d_2^{\kern 1pt 2}\,(x_3^2+x_1^2-d_2^{\kern 1pt 2})+d_3^{\kern 1pt 2}
\,(x_1^2+x_2^2-d_3^{\kern 1pt 2}),
\quad
\mytag{6.8}\\
\hskip -9em
\gathered
\tilde p_{\kern 1pt 8}=x_1^2\,d_1^{\kern 1pt 2}\,(x_2^2+x_3^2
-d_1^{\kern 1pt 2})+x_2^2\,d_2^{\kern 1pt 2}\,(x_3^2+x_1^2
-d_2^{\kern 1pt 2})\,+\\
\vphantom{1}\hskip 23em +\,x_3^2\,d_3^{\kern 1pt 2}\,(x_1^2
+x_2^2-d_3^{\kern 1pt 2}).
\endgathered
\mytag{6.9}
\endgather
$$
Using the formulas \mythetag{3.2}, \mythetag{6.3}, \mythetag{6.4}, 
\mythetag{6.5}, \mythetag{6.6}, \mythetag{6.7}, \mythetag{6.8}, 
and \mythetag{6.9}, 
now we can write the $S_3$ factor equations for the cuboid equations
\mythetag{1.1}. For this purpose it is convenient to use the
polynomials $p_1$, $p_2$, and $p_3$ from \mythetag{1.7}:
$$
\xalignat 2
&\hskip -2em
x_1^2+x_2^2+x_3^2-L^2=0,
&&p_{\kern 1pt 1}+p_{\kern 1pt 2}+p_{\kern 1pt 3}=0,
\quad\\
&\hskip -2em
d_1\,p_{\kern 1pt 1}+d_2\,p_{\kern 1pt 2}+d_3\,p_{\kern 1pt 3}=0,
&&x_1\,p_{\kern 1pt 1}+x_2\,p_{\kern 1pt 2}+x_3\,p_{\kern 1pt 3}=0,
\quad\\
\vspace{-0.7em}
\mytag{6.10}\\
\vspace{-0.7em}
&\hskip -2em
d_1^{\kern 1pt 2}\,p_{\kern 1pt 1}+d_2^{\kern 1pt 2}\,p_{\kern 1pt 2}
+d_3^{\kern 1pt 2}\,p_{\kern 1pt 3}=0,
&&x_1^2\,p_{\kern 1pt 1}+x_2^2\,p_{\kern 1pt 2}
+x_3^2\,p_{\kern 1pt 3}=0,
\quad\\
&\hskip -4em
\gathered
x_1\,d_1\,p_{\kern 1pt 1}+x_2\,d_2\,p_{\kern 1pt 2}\,+\\
\vphantom{1}\hskip 6em +\,x_3\,d_3\,p_{\kern 1pt 3}=0,
\endgathered
&&\hskip -2em
\gathered
x_1^2\,d_1^{\kern 1pt 2}\,p_{\kern 1pt 1}+x_2^2\,d_2^{\kern 1pt 2}
\,p_{\kern 1pt 2}\,+\\
\vphantom{1}\hskip 6em +\,x_3^3\,d_3^{\kern 1pt 2}\,p_{\kern 1pt 3}=0.
\endgathered
\quad
\endxalignat
$$ 
Since the polynomials $\tilde p_{\kern 1pt 1}$, $\tilde p_{\kern 1pt 2}$, 
$\tilde p_{\kern 1pt 3}$, $\tilde p_{\kern 1pt 4}$, $\tilde p_{\kern 1pt 5}$, 
$\tilde p_{\kern 1pt 6}$, $\tilde p_{\kern 1pt 7}$, $\tilde p_{\kern 1pt 8}$
constitute a basis of the ideal $I_{\text{sym}}$, the equations 
\mythetag{6.10} compose a complete set of $S_3$ factor equations.
\head
7. Comparison with the previously obtained factor equations.
\endhead
     In the previous paper \mycite{45} eight factor equations were already
derived. But they were written in so-called $E$-form, i\.\,e\. in terms of
the values of the elementary multisymmetric polynomials \mythetag{1.3},
\mythetag{1.4}, and \mythetag{1.5}. In order to compare the previously 
obtained equations from \mycite{45} with the equations \mythetag{6.10} we 
need to convert them into $xd$-form by means of the mapping $\varphi$
from \mythetag{2.1}.\par
     Let's consider the first of the previously obtained factor equations.
In its $E$-form this equation is written as follows (see \thetag{4.3} 
in \mycite{45}):
$$
\hskip -2em
E_{10}^2-2\,E_{20}-L^2=0. 
\mytag{7.1}
$$
In order to apply $\varphi$ to \mythetag{7.1} we should substitute 
$E_{10}=e_{\sssize [1,0]}$, $E_{20}=e_{\sssize [2,0]}$ and then use the
formulas \mythetag{1.3}. As a result we get the equation coinciding with 
the first equation in the left column of \mythetag{6.10}.\par
     The second of the previously obtained factor equations is the equation
\thetag{4.6} in \mycite{45}. In its $E$-form this equation is written as 
follows:
$$
\hskip -2em
E_{01}^2-2\,E_{02}-2\,L^2=0.
\mytag{7.2}
$$
Upon applying the mapping $\varphi$ to \mythetag{7.2} we get the equation
$$
\hskip -2em
2\,(x_1^2+x_2^2+x_3^2-L^2)-(p_{\kern 1pt 1}+p_{\kern 1pt 2}
+p_{\kern 1pt 3})=0,
\mytag{7.3}
$$
which is derived from the first equations in the left and right columns of 
\mythetag{6.10}.\par
     Let's proceed to the third of the previously obtained factor equations.
This is the equation \thetag{4.12} in \mycite{45}. In its $E$-form this equation 
is written as follows:
$$
\hskip -2em
2\,E_{12}+6\,E_{30}-2\,E_{01}\,E_{11}
+E_{10}\,E_{01}^2+3\,E_{10}\,L^2-E_{10}^3=0.
\mytag{7.4}
$$
Upon converting to the $xd$-form the equation \mythetag{7.4} looks like
$$
\hskip -2em
\gathered
e_{\sssize [1,0]}\,((p_{\kern 1pt 1}+p_{\kern 1pt 2}+p_{\kern 1pt 3})
-3\,(x_1^2+x_2^2+x_3^2-L^2))\,-\\
-\,2\,(x_1\,p_{\kern 1pt 1}+x_2\,p_{\kern 1pt 2}
+x_3\,p_{\kern 1pt 3})=0.
\endgathered
\mytag{7.5}
$$
It is easy to see that \mythetag{7.5} can be derived from the 
first and the second equations in the right column of \mythetag{6.10} and
from the first equation in the left column of \mythetag{6.10}.\par
     The fourth of the previously obtained factor equations is the equation
\thetag{4.19} in \mycite{45}. In its $E$-form this equation is written as 
follows:
$$
\hskip -2em
2\,E_{21}+6\,E_{03}-2\,E_{10}\,E_{11}
+E_{01}\,E_{10}^2+5\,E_{01}\,L^2-E_{01}^3=0.
\mytag{7.6}
$$
Upon applying the mapping $\varphi$ to \mythetag{7.6} we get the equation
$$
\hskip -2em
\gathered
e_{\sssize [0,1]}\,(3\,(p_{\kern 1pt 1}+p_{\kern 1pt 2}+p_{\kern 1pt 3})
-5\,(x_1^2+x_2^2+x_3^2-L^2))\,-\\
-\,2\,(d_1\,p_{\kern 1pt 1}+d_2\,p_{\kern 1pt 2}
+d_3\,p_{\kern 1pt 3})=0.
\endgathered
\mytag{7.7}
$$
The equation \mythetag{7.7} can be derived from the first and the second 
equations in the left column of \mythetag{6.10} and from the first equation 
in the right column of \mythetag{6.10}.\par
     The fifth of the previously obtained factor equations is more 
complicated. It is given by the formula \thetag{5.5} in \mycite{45}.
Here is its $E$-form: 
$$
\hskip -2em
\aligned
8\,E_{10}\,&E_{12}-8\,E_{01}\,E_{21}
-8\,E_{11}^2+4\,E_{01}^2\,E_{10}^2\,-\\
&-\,E_{01}^4-3\,E_{10}^4
+10\,E_{10}^2\,L^2+4\,E_{01}^2\,L^2+L^4=0.
\endaligned
\mytag{7.8}
$$
Upon converting to the $xd$-form the equation \mythetag{7.8} looks like
$$
\gathered
18\,(x_1^2\,p_{\kern 1pt 1}+x_2^2\,p_{\kern 1pt 2}+x_3^2\,p_{\kern 1pt 3})
+6\,(x_1\,d_1\,p_{\kern 1pt 1}+x_2\,d_2\,p_{\kern 1pt 2}+x_3\,d_3
\,p_{\kern 1pt 3})\,-\\
-\,8\,e_{\sssize [0,1]}\,(d_1\,p_{\kern 1pt 1}+d_2\,p_{\kern 1pt 2}
+d_3\,p_{\kern 1pt 3})-24\,e_{\sssize [1,0]}\,(x_1\,p_{\kern 1pt 1}
+x_2\,p_{\kern 1pt 2}+x_3\,p_{\kern 1pt 3})\,+\\
+\,(8\,e_{\sssize [2,0]}+3\,e_{\sssize [0,1]}^2+4\,e_{\sssize [1,0]}^2
+6\,e_{\sssize [0,2]})\,(p_{\kern 1pt 1}+p_{\kern 1pt 2}
+p_{\kern 1pt 3})\,+\\
+\,(2\,e_{\sssize [2,0]}-4\,e_{\sssize [0,1]}^2-11\,e_{\sssize [1,0]}^2
-L^2)\,(x_1^2+x_2^2+x_3^2-L^2)=0.
\endgathered
\quad
\mytag{7.9}
$$
Like \mythetag{7.5} and \mythetag{7.7}, the equation \mythetag{7.9} is
a linear combination of the equations \mythetag{6.10} with coefficients
in $\Sym\!\Bbb Q[M,L]$, i\.\,e\. it can be derived from \mythetag{6.10}. 
\par
     The next step is to consider the sixth of the previously obtained 
factor equations. It is given by the formula \thetag{5.10} in \mycite{45}.
Here is its $E$-form: 
$$
\hskip -2em
\aligned
-8\,E_{10}\,&E_{12}+8\,E_{01}\,E_{21}-8\,E_{11}^2+4\,E_{01}^2\,E_{10}^2\,-\\
&-\,E_{10}^4-3\,E_{01}^4+20\,E_{01}^2\,L^2-2\,E_{10}^2\,L^2-5\,L^4=0.
\endaligned
\mytag{7.10}
$$
The equation is similar to \mythetag{7.8} and is equally complicated as the
equation \mythetag{7.8} since it is of the same order with respect to its
variables. \pagebreak Upon converting to the $xd$-form the equation 
\mythetag{7.10} is written as follows: 
$$
\gathered
6\,(x_1^2\,p_{\kern 1pt 1}+x_2^2\,p_{\kern 1pt 2}+x_3^2\,p_{\kern 1pt 3})
+18\,(x_1\,d_1\,p_{\kern 1pt 1}+x_2\,d_2\,p_{\kern 1pt 2}+x_3\,d_3
\,p_{\kern 1pt 3})\,-\\
-\,24\,e_{\sssize [0,1]}\,(d_1\,p_{\kern 1pt 1}+d_2\,p_{\kern 1pt 2}
+d_3\,p_{\kern 1pt 3})-8\,e_{\sssize [1,0]}\,(x_1\,p_{\kern 1pt 1}
+x_2\,p_{\kern 1pt 2}+x_3\,p_{\kern 1pt 3})\,+\\
+\,(8\,e_{\sssize [2,0]}+9\,e_{\sssize [0,1]}^2-4\,e_{\sssize [1,0]}^2
+18\,e_{\sssize [0,2]})\,(p_{\kern 1pt 1}+p_{\kern 1pt 2}
+p_{\kern 1pt 3})\,-\\
-\,(10\,e_{\sssize [2,0]}20\,e_{\sssize [0,1]}^2-7\,e_{\sssize [1,0]}^2
-5\,L^2)\,(x_1^2+x_2^2+x_3^2-L^2)=0.
\endgathered
\quad
\mytag{7.11}
$$
Like \mythetag{7.9}, the equation \mythetag{7.11} is a linear combination 
of the equations \mythetag{6.10} with coefficients in $\Sym\!\Bbb Q[M,L]$, 
i\.\,e\. it can be derived from \mythetag{6.10}.\par
     Let's proceed to the seventh of the previously obtained factor equations.
This is the equation \thetag{5.17} in \mycite{45}. In its $E$-form this equation 
is written as follows:
$$
\hskip -2em
\gathered
4\,E_{11}\,E_{21}
-2\,E_{11}\,E_{01}^3
+6\,E_{12}\,E_{01}^2
+2\,E_{12}\,E_{10}^2
-\,E_{10}^3\,E_{01}^2\,+\\
+\,E_{10}\,E_{01}^4
-2\,E_{12}\,L^2
-E_{10}\,E_{01}^2\,L^2
+2\,E_{10}^3\,L^2
-2\,E_{10}\,L^4=0.
\endgathered
\quad
\mytag{7.12}
$$
Upon converting to the $xd$-form the equation \mythetag{7.12} looks like
$$
\gathered
-4\,e_{\sssize [1,0]}\,(d_1^{\kern 1pt 2}\,p_{\kern 1pt 1}
+d_2^{\kern 1pt 2}\,p_{\kern 1pt 2}+d_3^{\kern 1pt 2}\,p_{\kern 1pt 3})
-8\,e_{\sssize [0,1]}\,(x_1\,d_1\,p_{\kern 1pt 1}
+x_2\,d_2\,p_{\kern 1pt 2}\,+\\
+\,x_3\,d_3\,p_{\kern 1pt 3})
-(4\,e_{\sssize [2,0]}+2\,e_{\sssize [0,2]}-2\,e_{\sssize [1,0]}^2
-3\,e_{\sssize [0,1]}^2)\,(x_1\,p_{\kern 1pt 1}+x_2\,p_{\kern 1pt 2}\,+\\
+\,x_3\,p_{\kern 1pt 3})+8\,e_{\sssize [1,0]}\,e_{\sssize [0,1]}
\,(d_1\,p_{\kern 1pt 1}+d_2\,p_{\kern 1pt 2}+d_3\,p_{\kern 1pt 3})+
(3\,e_{\sssize [3,0]}\,-\\
-\,3\,e_{\sssize [1,1]}\,e_{\sssize [0,1]}
+3\,e_{\sssize [2,0]}\,e_{\sssize [1,0]}-3\,e_{\sssize [0,2]}
\,e_{\sssize [1,0]}-e_{\sssize [1,0]}^3
-e_{\sssize [1,2]}\,+\\
+\,L^2\,e_{\sssize [1,0]})\,(p_{\kern 1pt 1}+p_{\kern 1pt 2}
+p_{\kern 1pt 3})+(2\,e_{\sssize [1,2]}+2\,e_{\sssize [1,0]}^3
-8\,e_{\sssize [2,0]}\,e_{\sssize [1,0]}\,+\\
+\,2\,e_{\sssize [0,2]}\,e_{\sssize [1,0]}
+2\,e_{\sssize [1,0]}\,L^2)\,(x_1^2+x_2^2+x_3^2-L^2)=0.
\endgathered
\quad
\mytag{7.13}
$$
Again, looking at \mythetag{7.13}, we see that this equation is a
linear combination of the equations \mythetag{6.10} with coefficients 
in $\Sym\!\Bbb Q[M,L]$, i\.\,e\. it can be derived from 
\mythetag{6.10}.\par
     The eighth of the previously obtained factor equations is similar
to the seventh one. It is given by the formula \thetag{5.22} in \mycite{45}.
Here is its $E$-form: 
$$
\gathered
4\,E_{11}\,E_{12}-2\,E_{11}\,E_{10}^3+6\,E_{21}\,E_{10}^2
+2\,E_{21}\,E_{01}^2-E_{01}^3\,E_{10}^2+E_{01}\,E_{10}^4\,+\\
+\,2\,E_{21}\,L^2-2\,E_{11}\,E_{10}\,L^2
+2\,E_{01}\,E_{10}^2\,L^2+E_{01}^3\,L^2
-3\,E_{01}\,L^4=0.
\endgathered
\quad
\mytag{7.14}
$$
Converting the equation \mythetag{7.14} to the $xd$-form, we obtain 
$$
\gathered
-4\,e_{\sssize [0,1]}\,(x_1^2\,p_{\kern 1pt 1}+x_2^2\,p_{\kern 1pt 2}
+x_3^2\,p_{\kern 1pt 3})
-8\,e_{\sssize [1,0]}\,(x_1\,d_1\,p_{\kern 1pt 1}
+x_2\,d_2\,p_{\kern 1pt 2}\,+\\
+\,x_3\,d_3\,p_{\kern 1pt 3})+8\,e_{\sssize [1,0]}\,e_{\sssize [0,1]}
\,(x_1\,p_{\kern 1pt 1}+
x_2\,p_{\kern 1pt 2}+x_3\,p_{\kern 1pt 3})+(2\,e_{\sssize [1,0]}^2\,+\\
+\,2\,L^2)\,(d_1\,p_{\kern 1pt 1}+d_2\,p_{\kern 1pt 2}
+d_3\,p_{\kern 1pt 3})-(4\,e_{\sssize [1,1]}\,e_{\sssize [1,0]}
+4\,e_{\sssize [2,0]}\,e_{\sssize [0,1]}\,-\\
-\,3\,e_{\sssize [0,1]}\,e_{\sssize [1,0]}^2+3\,L^2\,e_{\sssize [0,1]}
-2\,e_{\sssize [2,1]})\,(p_{\kern 1pt 1}+p_{\kern 1pt 2}
+p_{\kern 1pt 3})-(4\,e_{\sssize [2,1]}\,+\\
+\,6\,e_{\sssize [0,3]}-4\,e_{\sssize [2,0]}\,e_{\sssize [0,1]}
-4\,e_{\sssize [1,1]}\,e_{\sssize [1,0]}-3\,L^2\,e_{\sssize [0,1]}\,+\\
+\,5\,e_{\sssize [0,1]}\,e_{\sssize [1,0]}^2)\,(x_1^2+x_2^2+x_3^2-L^2)=0.
\endgathered
\quad
\mytag{7.15}
$$
Like \mythetag{7.13}, the equation \mythetag{7.15} is a linear combination 
of the equations \mythetag{6.10} with coefficients in $\Sym\!\Bbb Q[M,L]$.
This means that it can be derived from \mythetag{6.10}.\par
\head 
8. Concluding remarks.
\endhead
     The theorem~\mythetheorem{6.1} is the main result of this paper. 
It yields a basis for the ideal $I_{\text{sym}}$ and a complete list 
\mythetag{6.10} of the cuboid factor equations in $xd$-form. As we noted 
above the equation \mythetag{7.1} is equivalent to the first equation 
\mythetag{6.10}. Looking attentively at \mythetag{7.3}, \mythetag{7.5}, 
\mythetag{7.7}, \mythetag{7.9}, \mythetag{7.11}, \mythetag{7.13}, and 
\mythetag{7.15}, we find that the first seven equations in \mythetag{6.10} 
can be derived from the previously obtained eight factor equations 
\mythetag{7.1}, \mythetag{7.2}, \mythetag{7.4}, \mythetag{7.6}, 
\mythetag{7.8}, \mythetag{7.10}, \mythetag{7.12}, and \mythetag{7.14}. 
The last equation \mythetag{6.10} is new. Upon converting to an $E$-form 
it can be added to the list of previously obtained factor equations. 
However, this is not enough for to complete the list. The matter is that 
in $E$-form a complete list should include kernel equations. Therefore 
the equations $q_{\kern 1pt 1}=0$, $q_{\kern 1pt 2}=0$, $q_{\kern 1pt 3}=0$, 
$q_{\kern 1pt 4}=0$, $q_{\kern 1pt 5}=0$, $q_{\kern 1pt 6}=0$, 
$q_{\kern 1pt 7}=0$ given by the kernel polynomials \mythetag{2.4}, 
\mythetag{2.5}, \mythetag{2.6}, \mythetag{2.7}, \mythetag{2.8}, \mythetag{2.9}, 
and \mythetag{2.10} should be added.\par
\head 
9. Acknowledgments.
\endhead
    I am grateful to my colleague I.~Yu.~Cherdantsev who recommended me the
book \mycite{69} for learning Gr\"obner bases and their applications.

\enddocument
\end